\documentclass[12pt]{article}

\usepackage{epsf}
\usepackage{psfrag}

\usepackage{amsmath}
\usepackage{amsthm}
\usepackage{amssymb}

\theoremstyle{plain}
\newtheorem{thm}{Theorem}[section]
\newtheorem{prop}[thm]{Proposition}
\newtheorem{cor}[thm]{Corollary}

\newtheorem{intlemnp}[thm]{Lemma}

\newtheorem{intthmnp}[thm]{Theorem}
\newenvironment{thmnp}{\begin{intthmnp}}{\qed \end{intthmnp}}

\theoremstyle{definition}

\newtheorem{ex}[thm]{Example}

\theoremstyle{remark}
\newtheorem{rem}[thm]{Remark}
\newtheorem{prob}[thm]{Problem}

\renewcommand{\epsfsize}[2]{0.8\textwidth}

\def\thmref#1{Theorem~\ref{#1}}
\def\exref#1{Example~\ref{#1}}

\def\propref#1{Proposition~\ref{#1}}
\def\propertyref#1{Property~\ref{#1}}

\def\figref#1{Figure~\ref{#1}}

\def\secref#1{Section~\ref{#1}}
\def\appref#1{Appendix~\ref{#1}}
\def\propertyref#1{Property~\ref{#1}}
\def\presref#1{Presentation~(\ref{#1})}

\def\Mod{\text{Mod}}

\def\jmt{{\tt jmt}}
\def\jsnap{{\tt jsnap}}
\def\snappea{{\tt Snap\-Pea}}

\def\co{\colon\thinspace}

\title{Computing Triangulations of Mapping Tori of Surface Homeomorphisms}
\date{}
\author{Peter Brinkmann\footnote{%
This research was partially conducted by the
first author for the Clay Mathematics Institute.}
~and Saul Schleimer}

\begin{document}
\renewcommand{\thefootnote}{\null}
\maketitle%
\setcounter{footnote}{0}
\renewcommand{\thefootnote}{\arabic{footnote}}

\begin{abstract}
We present the mathematical background of a software package that
computes triangulations of mapping tori of surface homeomorphisms, suitable
for Jeff Weeks's program \snappea. The package is an extension of the
software described in \cite{pbexp}. It consists of two programs. \jmt~
computes triangulations and prints them in a human-readable format.
\jsnap~ converts this format into \snappea's triangulation file
format and may be of independent interest because it allows for quick
and easy generation of input for \snappea.
As an application, we obtain a new solution to the restricted conjugacy
problem in the mapping class group.
\end{abstract}

\section{Introduction} \label{intro}
In \cite{pbexp}, the first author described a software
package that provides an
environment for computer experiments with automorphisms of surfaces with one
puncture.
The purpose of this paper is to present the mathematical background of an
extension of this package that
computes triangulations of mapping tori of such homeomorphisms, suitable
for further analysis with Jeff Weeks's program \snappea~\cite{snappea}.\footnote{Software
available at {\tt http://thames.northnet.org/weeks/index/SnapPea.html}}

Pseudo-Anosov homeomorphisms are of particular interest because their mapping
tori are hyperbolic $3$-manifolds of finite volume \cite{thurston2}.
The software described in \cite{pbexp} recognizes pseudo-Anosov homeomorphisms.
Combining this with the programs discussed here, we obtain a powerful tool
for generating and analyzing large numbers of hyperbolic $3$-manifolds.

The software package described in \cite{pbexp} takes an automorphism $\phi$
of a surface $S$ with one puncture (given as a sequence of Dehn twists) and
computes the induced outer automorphism of the fundamental group of $S$,
represented by a homotopy equivalence $f\co G\rightarrow G$ of a finite
graph $G\subset S$ homotopy equivalent to $S$, together with a
loop $\sigma$ in $G$ homotopic to a loop around the puncture of $S$.
The map $f$ and the loop $\sigma$ determine $\phi$ up
to isotopy \cite[Section~5.1]{pbexp}.

In \secref{triang}, we describe an effective algorithm for computing
a triangulation of the mapping torus of $\phi\co S\rightarrow S$,
given only $f\co G\rightarrow G$ and $\sigma$ (\thmref{main}).
We also present an
analysis of the complexity of this algorithm (\propref{complexity}).
The first part of the
software package is a program (called \jmt) that implements
this procedure. The program \jmt~prints its output in an intermediate
human-readable format.

In \secref{conjref}, we explain how to use the software discussed here and
the isometry checker of \snappea~to solve the restricted conjugacy problem
in the mapping class group (i.e., the question of whether two pseudo-Anosov
homeomorphisms are conjugate in the mapping class group). This problem
was previously solved in \cite{MR89f:57016} and \cite{MR80f:57003}. One
distinguishing feature of our solution is that much of it has already
been implemented.

\appref{format} discusses the second program in the software package (called
\jsnap), which converts the intermediate format of \jmt~into
\snappea's triangulation file format. Since \snappea's format is rather
complicated, it is not easy to generate input files for \snappea, and
\jsnap~may be of independent interest
because it allows users to generate input for \snappea~without having
to understand \snappea's file format.

Finally, in \appref{exsec}, we present some sample computations that exhibit
some of the capabilities of the combination of \snappea~and the software
discussed here.

Immediate applications of the software described here include an
experimental investigation of possible relationships between dynamical
properties of pseudo-Anosov homeomorphisms (as computed by the first
author's train track software) and topological properties of their
mapping tori (as computed by \snappea). For example, one might look for
a relationship between growth rate and volume.  Another area where the
package described in this paper has already been used is the study
of slalom knots as introduced by Norbert A'Campo \cite{acampo1}.

The software package is written in Java and should be universally portable.
The programs \jmt~and \jsnap~are command line software and can be used to
examine a large number of examples as a batch job. A graphical user interface
with an online help feature is also available.

The package, including binary files, source code, complete online
documentation, and a user manual, is available at
{{\tt http://www.math.uiuc.edu/ \~{}brinkman/}}.

We would like to thank Mladen Bestvina and John Stallings 
for many helpful discussions,
as well as Jeff Weeks and Bill Floyd for explaining \snappea's
intricacies. We would also like to express our gratitude to Kai-Uwe Bux
for critiquing an early version of this paper.

\section{Computing triangulations}\label{triang}

Let $\phi\co S\rightarrow S$ be an automorphism of a surface $S$ with one
puncture, represented by a homotopy equivalence $f\co G\rightarrow G$ of
a finite graph $G$ and a loop $\sigma$ in $G$ representing a loop
around the puncture of $S$ (see \secref{intro}).
There is no loss in assuming that $f\co G\rightarrow G$ maps vertices to
vertices and that the restriction of $f$ to the interior of each edge
of $G$ is an immersion.
 
In this section, we outline an effective procedure that computes a 
triangulation of the mapping torus of $\phi$ given only $f$ and $\sigma$.
To this end, we construct a simplicial $2$-complex $K$ and a face
pairing $e$ with the following properties.
\begin{enumerate}
\item \label{prop1} The space $|K|$ is homeomorphic to a torus.
\item For each $2$-simplex $\Delta$ of $K$, there exists a $2$-simplex
$\Delta'$ of $K$ and an orientation reversing simplicial homeomorphism
$e^{\null}_\Delta\co \Delta\rightarrow \Delta'$ such that
$e^{-1}_{\Delta'}=e^{\null}_\Delta$.
\item The space $K/e$ is homotopy equivalent to the mapping torus of $f$.


\item \label{prop2} If we let $M=(\text{cone over~} K)/e$ and obtain $M'$
from $M$ by removing the cone point, then $M'$ is a $3$-manifold (in
particular, the links of vertices in $M'$ are 2-spheres).\label{3mfcond}
\end{enumerate}

In this situation, $M'$ is homotopy equivalent to the mapping torus of $f$,
which in turn is homotopy equivalent to the mapping torus $M_\phi$ of $\phi$.
As $M'$ is a $3$-manifold, $M'$ is homeomorphic to $M_\phi$
\cite[page~6]{johann}.

The triangulation of $K$ induces a triangulation of $M$, i.e., the tetrahedra
of $M$ are cones over the triangles of $K$. The vertices of $K$ give rise to
finite vertices of $M$, and the cone point is an ideal vertex corresponding
to the torus cusp of $M_\phi$. By computing the links of vertices,
\snappea~recognizes finite vertices (whose links are 2-spheres) and ideal
vertices (whose links are tori or Klein bottles).

Hence, we have reduced to problem of constructing a triangulation of the
mapping torus of $\phi$ to the construction of the 2-complex $K$ and face
pairing $e$, given only the homotopy equivalence $f:G\rightarrow G$ and
the loop $\sigma$.
The construction of $K$ and $e$ is the purpose of the remainder of this
section.


The construction of $K$ and $e$ proceeds in two steps. We construct the
$2$-torus $T$ by gluing annuli using 
Stallings's folding construction \cite{stallings}.
Then we construct a triangulation and a face pairing for each of the annuli.

\subsection{Step 1: Subdividing and folding}\label{foldsubsec}

We review the notion of subdividing and folding \cite{stallings,hb1}.
Let $G,G'$ be finite graphs, and let $f\co G'\rightarrow G$ be a map that maps
vertices to vertices and edges to edge paths.

If $f$ fails to be an
immersion, then there exist two distinct edges $a,~b$ in $G'$ emanating
from the same vertex such that
$f(a)$ and $f(b)$ have a nontrivial initial path in common. 
We construct a new graph $G_1'$ by subdividing $a$ (resp.\ $b$) into two
edges $a_1,~a_2$ (resp.\ $b_1,~b_2$).

\begin{figure}[tb]
\renewcommand{\epsfsize}[2]{\textwidth}
\psfrag{a}{$a$}
\psfrag{b}{$b$}
\psfrag{c}{$c$}
\psfrag{d}{$d$}
\psfrag{a1}{$a_1$}
\psfrag{b1}{$b_1$}
\psfrag{a2}{$a_2$}
\psfrag{b2}{$b_2$}
\psfrag{a1eqb1}{$a_1=b_1$}
\psfrag{subdivision}{subdivision}
\psfrag{fold}{fold}
\centerline{\epsfbox{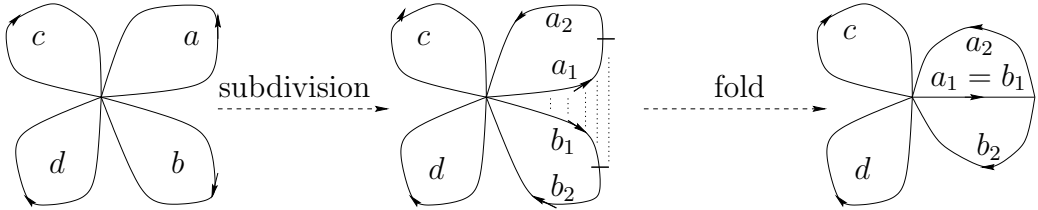}}
\caption{Subdividing and folding.}\label{stfoldpic}
\end{figure}

Now $f$ factors through $G_1'$, i.e., there are maps
$s\co G'\rightarrow G_1'$ and $g\co G_1'\rightarrow G$ such that
$f=g\circ s$. Moreover, we can choose $s$ and $g$ such that $s(a)=a_1a_2$,
$s(b)=b_1b_2$ and $g(a_1)=g(b_1)$. We obtain a new graph $G_2'$ from $G_1'$
by identifying the edges $a_1$ and $b_1$. Then $g$ factors through
$G_2'$, i.e., there
is a map $h\co G_2'\rightarrow G$ such that $g=h\circ p$, where $p$ is the
natural projection $p\co G_1'\rightarrow G_2'$ (see \figref{stfoldpic}). We
refer to this process as {\em folding} $a_1$ and $b_1$.

\begin{rem}
The notion of folds used in \cite{hb1} differs slightly from that
introduced in \cite{stallings}. In \cite{hb1}, the authors consider
homotopy equivalences $f:G\rightarrow G$, and folding changes {\em both}
the domain and the range of $f$, whereas in \cite{stallings}, the author
considers maps $f:G'\rightarrow G$, and folding only affects the domain $G'$.
The notion of folds used in this paper is a slight modification of the
folds in \cite{stallings}.
\end{rem}

The remainder of this subsection details how to
construct a sequence of graphs and maps
$$ G=G_0 \overset{s_0}{\longrightarrow} G_1 \overset{p_0}{\longrightarrow} G_2
\ldots G_{2n-2} \overset{s_{n-1}}{\longrightarrow} G_{2n-1}
\overset{p_{n-1}}{\longrightarrow} G_{2n} \overset{g_{n}}{\longrightarrow} G$$
such that
$$f=g_{n}\circ p_{n-1} \circ s_{n-1} \circ \ldots \circ p_0 \circ s_0,$$
where $s_i\co G_{2i}\rightarrow G_{2i+1}$ is a subdivision,
$p_i\co G_{2i+1}\rightarrow G_{2i+2}$ is a Stallings fold,
and $g_{n}\co G_{2n}\rightarrow G_0$ is an immersion.
Since $f$ is a homotopy equivalence, $g_{n}$ will be onto, hence a
homeomorphism. Moreover, for each $i=0,\ldots,2n$, we will construct a
loop $\sigma_i$ in $G_i$ corresponding to a loop around the puncture of $S$.

Let
$f\co G\rightarrow G$ be induced by a homeomorphism $\phi\co S\rightarrow S$,
and let $\sigma$ denote an edge loop in $G$ corresponding to a loop around the
puncture of $S$. Let $G=G_0$, $g_0=f\co G_0\rightarrow G$, and
$\sigma_0=\sigma$.

Suppose that $g_0$ is not an immersion. Then there exist two edges $a,~b$
emanating from the same vertex in $G_0$ such that $g_0(a)$ and $g_0(b)$ have
a common initial segment. Since $g_0$ is induced by the homeomorphism
$\phi\co S\rightarrow S$, we can find $a$ and $b$ such that $a$ and $b$ are
adjacent in the embedding of $G_0$ in $S$.

Since the loop $\sigma_0$ in $G$ is homotopic to a loop around the puncture,
$a$ and $b$ will be adjacent in the spelling of $\sigma_0$. Hence, we can
detect $a$ and $b$ algorithmically by looking for cancellation between the
images of adjacent edges in the spelling of $\sigma_0$.

We obtain $G_1$ from $G_0$ by subdividing
$a$ and $b$, and we obtain $G_2$ from $G_1$ by folding the initial segments
of $a$ and $b$.
As above, we construct maps $s_0\co G_0\rightarrow G_1$,
$p_0\co G_1\rightarrow G_2$,
and $g_1\co G_2\rightarrow G$ such that $g_0=g_1\circ p_0\circ s_0$.
Let $\sigma_1 = s_0(\sigma_0)$ and obtain $\sigma_2$ from $p_0(\sigma_1)$
by tightening.
Since the edges $a$ and $b$
are adjacent in the embedding of $G$ in $S$, the embedding of $G_0$ in $S$
induces an embedding of $G_1$ and $G_2$ in $S$, and $\sigma_1$ and $\sigma_2$
are homotopic to $\sigma_0$ in $S$.

Note that the {\em size} of $g_1$, i.e., the sum of the lengths of the
images under $g_1$ of the edges in $G_2$, is strictly smaller than the
size of $g_0$.
Hence, after repeating this construction finitely many times, we reach
a map $g_{n}\co G_{2n}\rightarrow G$ that cannot be folded and thus has
to be an immersion. We have found the desired sequence of subdivisions
and folds.

\begin{ex}\label{ex1}
Let $G$ be the graph
with one vertex and four edges, labeled $a,\ldots,d$, embedded in a punctured
surface $S$ of genus $2$ as shown in \figref{embed}. We consider the
following homotopy equivalence $f\co G\rightarrow G$ induced by an
automorphism of $S$.
\begin{figure}[tb]
\renewcommand{\epsfsize}[2]{0.45\textwidth}
\psfrag{a}{$a$}
\psfrag{b}{$b$}
\psfrag{c}{$c$}
\psfrag{d}{$d$}
\centerline{\epsfbox{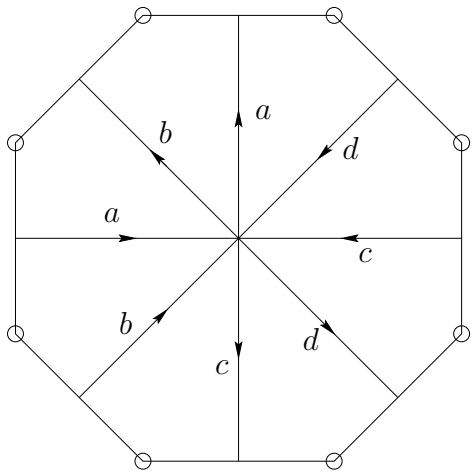}}
\caption{The graph $G$ embedded in the surface $S$. The corners of the octagon
correspond to the puncture of $S$. Two faces of the octagon are glued via an
orientation-reversing map if the edge labels match up. }
\label{embed}
\end{figure}

\begin{eqnarray*}
f(a)&=&a\underline{cd\bar c\bar b}\\
f(b)&=&bc\bar d\bar cb\underline{cd\bar c\bar b}\\
f(c)&=&bc\bar d\bar d\\
f(d)&=&dd\bar c\bar bd\\
\sigma&=&a\bar b\bar abc\bar d\bar cd\\
\end{eqnarray*}
In $f(\sigma)$, cancellation occurs between the underlined parts of $f(a)$
and $f(\bar b)$, and we subdivide $a$ and $b$ in preparation for folding, which
gives us the maps $s_0\co G\rightarrow G_1$ and $g_1\co G_1\rightarrow G$
(see \figref{stfoldpic}).
In order to reduce notational complexity, we only change the labels of
those edges that are subdivided.
\begin{eqnarray*}
s_0(a)&=&a_1a_2\\
s_0(b)&=&b_1b_2\\
s_0(c)&=&c\\
s_0(d)&=&d\\
\sigma_1&=&a_1a_2\bar b_2\bar b_1\bar a_2\bar a_1b_1b_2c\bar d\bar cd\\
\end{eqnarray*}

Now we fold the edges $\bar{a_2}$ and $\bar{b_2}$.
\begin{eqnarray*}
p_0(a_1)&=& a_1\\
p_0(a_2)&=& b_2\\
p_0(b_1)&=& b_1\\
p_0(b_2)&=& b_2\\
p_0(c)&=&c\\
p_0(d)&=&d\\
\end{eqnarray*}

Finally, we compute the map $g_1\co G_2\rightarrow G$.

\begin{eqnarray*}
g_1(a_1)&=& a\\
g_1(b_1)&=& bc\bar d\bar cb\\
g_1(b_2)&=& cd\bar c\bar b\\
g_1(c)&=& bc\bar d\bar d\\
g_1(d)&=& dd\bar c\bar bd\\
\sigma_2&=&a_1\bar b_1\bar b_2\bar a_1b_1b_2c\bar d\bar cd\\
\end{eqnarray*}

\end{ex}

\subsection{Step 2: Triangulating annuli}\label{triangsec}

We can interpret the loops $\sigma_i$ as immersions
$\sigma_i\co S^1\rightarrow G_i$. The preimage of the vertex set of $G_i$
subdivides $S^1$ into intervals, and the restriction of $\sigma_i$ to such
an interval is a homeomorphism onto the interior of an edge in $G_i$. Hence,
we can label each interval with the corresponding edge in $G_i$.
We refer to this construction as {\em spelling $\sigma_i$
along $S^1$}.

Now, for each $i\in\{0,\ldots,2n\},$ we take an annulus $A_i$ and
spell the word $\sigma_i$ along
one boundary component and $\sigma_{i+1}$ along the other.
We orient the two boundary components of $A_i$ such that they are
freely homotopic as oriented loops.

\begin{figure}[tb]
\renewcommand{\epsfsize}[2]{0.5\textwidth}
\psfrag{fold}{fold}
\psfrag{subdivision}{subdivision}
\psfrag{finalgl}{final gluing}
\centerline{\epsfbox{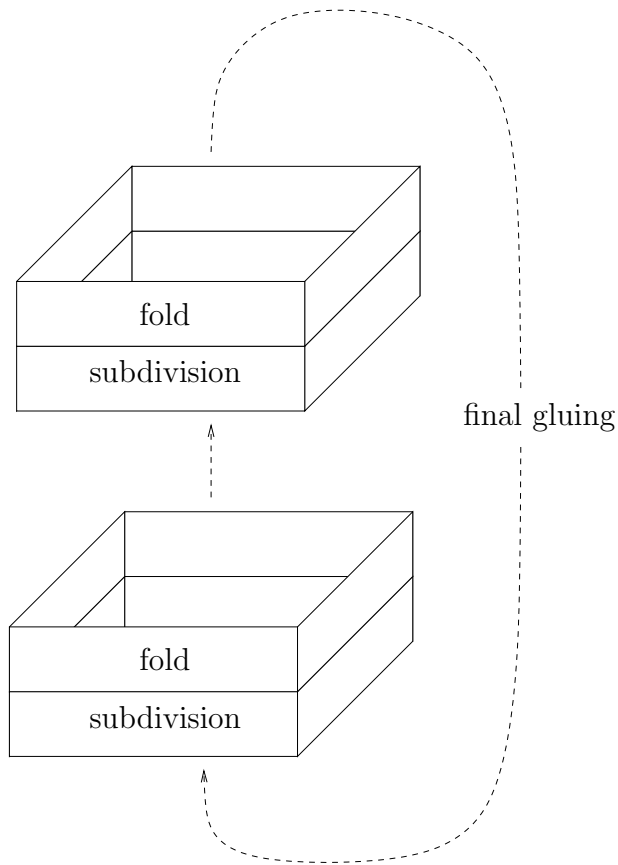}}
\caption{Decomposition into annuli.}\label{layerpic}
\end{figure}

\begin{figure}[tb]
\renewcommand{\epsfsize}[2]{0.9\textwidth}
\psfrag{A0}{$A_0$}
\psfrag{A1}{$A_1$}
\psfrag{a}{$a$}
\psfrag{b}{$b$}
\psfrag{c}{$c$}
\psfrag{d}{$d$}
\psfrag{a1}{$a_1$}
\psfrag{a2}{$a_2$}
\psfrag{b1}{$b_1$}
\psfrag{b2}{$b_2$}
\psfrag{D0}{$\Delta_0$}
\psfrag{D0p}{$\Delta_0'$}
\psfrag{D1}{$\Delta_1$}
\psfrag{D1p}{$\Delta_1'$}
\psfrag{D2}{$\Delta_2$}
\psfrag{D2p}{$\Delta_2'$}
\centerline{\epsfbox{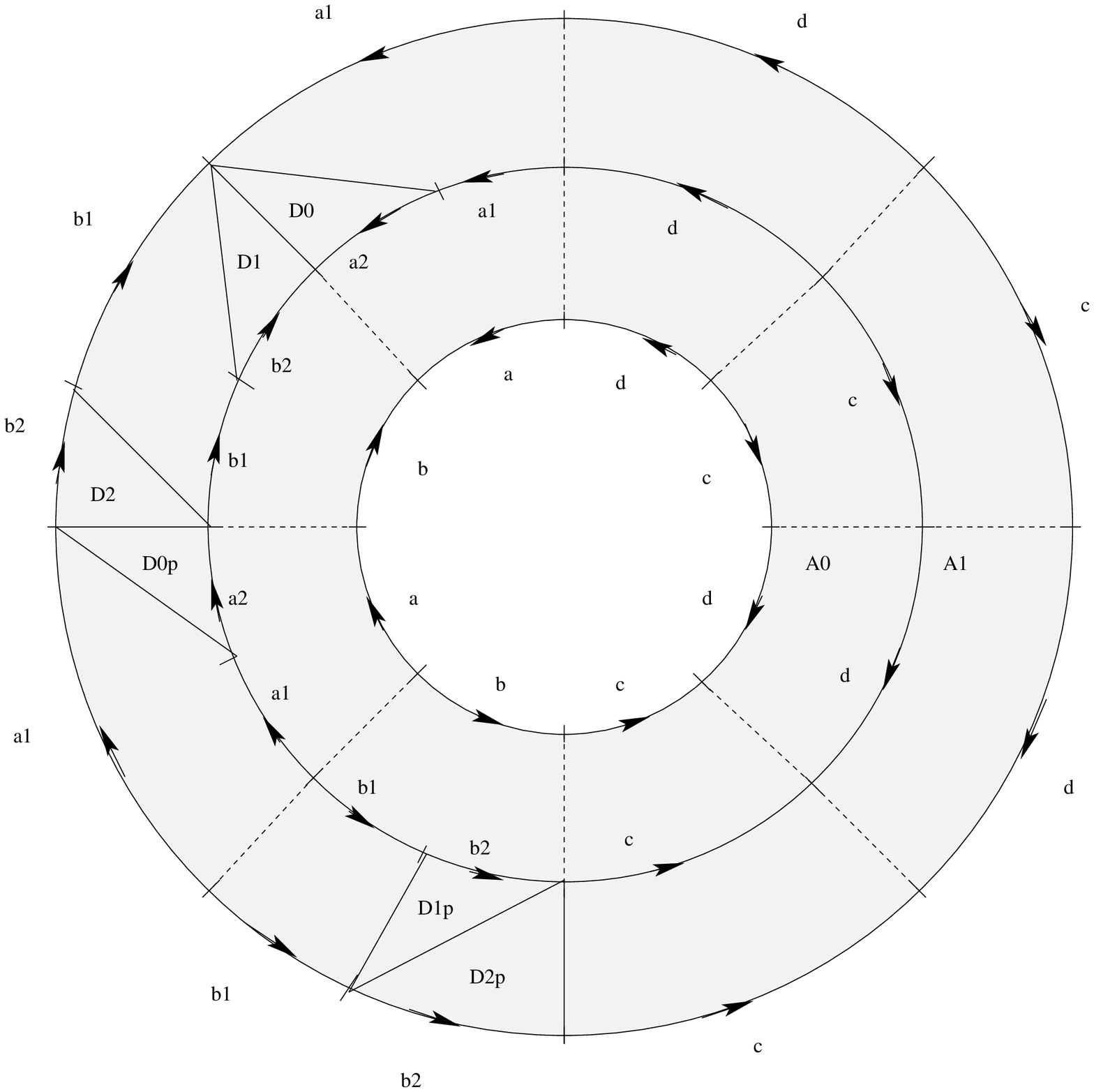}}
\caption{The gluing of $A_0$ and $A_1$ for \exref{ex1}.}\label{exglfig}
\end{figure}

This labeling
defines a gluing of $A_i$ and $A_{i+1}$, and the homeomorphism
$g_n\co G_{2n}\rightarrow G=G_0$ induces a gluing of $A_{2n-1}$ and $A_0$
(which we refer to as the {\em final gluing}),
giving us the desired torus $T$ (see \figref{layerpic}). \figref{exglfig}
shows the gluing of $A_0$ and $A_1$ for \exref{ex1}.

Fix some $i\in\{0,\ldots,2n\}$ and spell $\sigma_i$ along $S^1$. Notice that
each label on $S^1$ occurs twice, once for each direction (see \exref{ex1} and
\figref{embed}). If we identify
corresponding intervals, we obtain the graph $G_i$, and if we take the cone
over $S^1$, remove the cone point, and identify intervals with identical
labels, then we obtain the surface $S$ \cite[Section~5.1]{pbexp}.

Hence, we only need to extend the edge pairing on the boundary of the
annulus $A_i$ to an appropriate face pairing of a triangulation of all of
$A_i$. Recall that we want to choose the face pairing in such a way that
the cone over $T$ (with the cone point removed) becomes a $3$-manifold
when we identify corresponding triangles.

We first find a suitable triangulation of an annulus $A_{2i}$
corresponding to a subdivision $s_i$. We decompose $A_{2i}$ into rectangles
corresponding to edges that are not subdivided and pentagons corresponding
to those edges that are subdivided
(see \figref{exglfig}). The edge pairing on the boundary of $A_{2i}$ induces
a pairing of rectangles (resp.\ pentagons).
Finding a triangulation of $A_{2i}$ that is compatible with this pairing is
straightforward.

We now construct a suitable triangulation of an annulus $A_{2i+1}$
corresponding to a fold $p_i$. Edges that are not involved in the fold give
rise to paired rectangles contained in $A_{2i+1}$ (see \figref{exglfig}),
and as above, we easily find a triangulation of these rectangles that is
compatible with the pairing.

Let $a, b$ denote the two edges involved in the fold, i.e., we have
$p_i(a)=p_i(b)=b'$. By exchanging $a$ and $b$ or reversing the orientation
of $a$ and $b$ as necessary, we may assume that the loop $\sigma_{2i+1}$
has a subpath of the form $w=a\bar{b}u\bar{a}$ or $w=a\bar{b}ub$, where $u$
is a path that contains neither $a$ nor $b$. For concreteness, we focus on
the case $w=a\bar{b}u\bar{a}$. The first fold of \exref{ex1} falls into
this case, with $w=a_2\bar{b_2}\bar{b_1}\bar{a_2}$.
The construction in the remaining case is similar.

The loop $\sigma_{2i+2}$ has a corresponding subpath of the form
$w'=u'\bar{b'}$.  Let $\Delta_0$ be the triangle spanned by the initial
endpoint of $w'$ and the occurrence of $a$ in $w$. We pair $\Delta_0$ with
the triangle $\Delta_0'$ spanned by the terminal endpoint of $w'$ and
the occurrence of $\bar{a}$ in $w$ (see \figref{exglfig}).

Let $\Delta_1$ be the triangle spanned by the occurrence of $\bar{b}$ in
$w$ and the initial endpoint of $w'$, and let $\Delta_2$ be the triangle
spanned by the occurrence of $\bar{b'}$ in $w'$ and the terminal endpoint
of $u$. Observe that after identifying paired edges, $\Delta_1$ and $\Delta_2$
have a side in common, so we can think of $\Delta_1$ and $\Delta_2$ as spanning
a rectangle between $\bar{b}$ and $\bar{b'}$, which induces a triangulation of
the rectangle spanned by the occurrence of $b$ in $\sigma_{2i+1}$ and the 
occurrence of $\bar{b}$ in $\sigma_{2i+2}$ (see \figref{exglfig}).
This completes the triangulation
and face pairing of $A_{2i+1}$, which completes our construction.

\begin{figure}[tb]
\renewcommand{\epsfsize}[2]{0.9\textwidth}
\psfrag{a1}{$a_1$}
\psfrag{a2}{$a_2$}
\psfrag{a1}{$a_1$}
\psfrag{a1}{$a_1$}
\psfrag{b1}{$b_1$}
\psfrag{b2}{$b_2$}
\psfrag{b1}{$b_1$}
\psfrag{b2}{$b_2$}
\psfrag{b1}{$b_1$}
\psfrag{b2}{$b_2$}
\centerline{\epsfbox{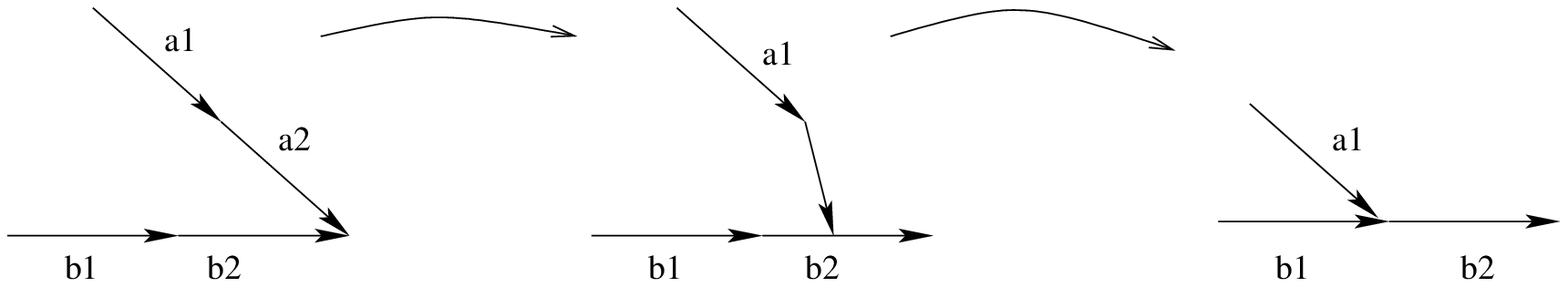}}
\caption{Continuous fold for \exref{ex1}.}
\label{contfold}
\end{figure}

Given the triangulation of $A_{2i+1}$ constructed above, we can think of
the annulus $A_{2i+1}$ as interpolating between the graph $G_{2i+1}$ and
the graph $G_{2i+2}$. In other words, we think of the fold as occurring
continuously, by identifying larger and larger segments of the edges $a$
and $b$ (see \figref{contfold}). Moreover, this continuous folding process
is compatible with the embedding of the graphs in the surface $S$. This
observation shows that the complex $K$ and its triangulation have the
desired properties, in particular \propertyref{prop2}.

Summing up, we have obtained the following main result of this paper.

\begin{thmnp}\label{main}
Let $f$ be a homotopy equivalence of a finite graph $G$, representing
a homeomorphism $\phi$ of a once-punctured surface. Then the following
is an effective procedure for computing a triangulation of the mapping
torus of $\phi$.

\begin{enumerate}
\item Decompose the homotopy equivalence $f\co G\rightarrow G$ into a sequence
of subdivisions and folds, followed by a homeomorphism (\secref{foldsubsec}).
\item Obtain the torus $K$ as a gluing of one annulus for each subdivision
and fold in the above decomposition (\figref{layerpic}).
\item Triangulate the individual annuli and construct a face pairing
(\secref{triangsec}).
\item Construct a triangulation of the mapping torus of the surface
homeomorphism $\phi$ by taking the cone of $K$ and glue tetrahedra
according to the face pairing.
\end{enumerate}

\end{thmnp}

The program \jmt~is an implementation of this procedure.

\subsection{Complexity analysis}

The purpose of this subsection is to obtain an estimate on the number
of tetrahedra in the triangulation that we have constructed.

Recall that we defined the size of a map $h\co G'\rightarrow G$ to be
the sum of the lengths (in the usual path metric) of the images
of the edges of $G'$. For example, the size of the first
map in \exref{ex1} is 23. Let $S(h)$ denote the size of $h$.

We say that a map $g\co G'\rightarrow G$ is {\em tight} if
\begin{enumerate}
\item for every edge $e$ of $G'$, the restriction of $g$ to the
interior of $e$ is an immersion, and
\item for every vertex $v$ of $G'$, there are two edges emanating from
$v$ that cannot be folded (not even after a subdivision).
\end{enumerate}
Note that tightness can always be achieved by homotopy.

\begin{prop}\label{complexity}
Let $f:G\rightarrow G$ be a tight homotopy equivalence representing a
homeomorphism $\phi$ of a once-punctured surface of genus $g$, and
assume that $G$ has no vertices of valence less than three.
Then the number of tetrahedra in the triangulation
of $M_\phi$ is bounded by $16(5g-2)S(f)$.
\end{prop}

\begin{proof}

Since folding reduces size and annuli come in pairs (a subdivision
annulus followed by a folding annulus), the number of annuli is bounded
above by $2S(f)$.  A simple application of Euler characteristics shows
that $G$ has no more than $6g-3$ edges because the valence of each vertex
is at least three. Subdividing and folding, however, may create additional
edges as well as vertices of valence one or two, so we need to understand
the effect of subdividing and folding on the number of edges.

To this end, we introduce the notion of {\em partial folds}, i.e.,
folds where both participating edges have to be subdivided, and
{\em full folds}, i.e., folds where at least one of the participating
edges is not subdivided \cite{hb1}. Clearly, a subdivision followed
by a full fold does not increase the number of edges, whereas a
subdivision followed by a partial fold increases the number of edges
by one.

A partial fold reduces the number of possible folds by one
because the map resulting from it is an immersion around
the new vertex created by the fold, so the only folds that
are possible after a partial fold are those that were available before.

Similarly, a full fold does not increase the number of possible folds.
This means that the number of partial folds that occurs in our construction
is bounded by the number of folds that the map $f\co G\rightarrow G$
admits. Since $f$ is tight, the number of folds at one vertex $v$ is bounded
by $\text{val}(v)-2$. Summing up, we see that the number of possible folds
is bounded by 
\[
	\sum_{v\in G^{(0)}}(\text{val}(v)-2)=-2\chi(G)=4g-2.
\]

Hence, the number of edges after a fold is bounded by $6g-3+4g-2=10g-5$.
Subdivisions increase the number of edges, so the number of edges at any
point in our construction is bounded by $10g-4.$

The number of tetrahedra belonging to one annulus
is bounded above by $8(5g-2)$ (four tetrahedra per edge), which gives
us a theoretical upper bound of $16(5g-2)S(f)$ on the number of tetrahedra
in our triangulation of $M_\phi$.
\end{proof}

Similar arguments show that for fixed genus, the time it takes to compute
a triangulation is linear in the size of the input.  We note that in
practice, partial folds seldom occur, and triangulations
tend to be considerably smaller than the bound given in \propref{complexity}.


\section{Solving the conjugacy problem}\label{conjref}

The algorithm from \secref{triang} and \snappea's isometry checker
provide a practical way of testing a necessary condition for two
pseudo-Anosov homeomorphisms $\phi, \gamma$ to be conjugate in the mapping
class group. Namely, we can compute the mapping tori of $\phi$ and $\gamma$
as in \secref{triang}, and then \snappea's isometry checker will determine
whether the two mapping tori are isometric. If they are not isometric,
we can immediately conclude that the $\phi$ and $\gamma$ are not conjugate.
However, if the mapping tori are isometric, we cannot yet conclude that
they are conjugate. The purpose of this section is to provide an effective
sufficient criterion for conjugacy.

\begin{prob}
Let $\gamma$ and $\phi$ be automorphisms of a surface 
$S$ with one puncture. Both automorphisms are assumed to be 
presented as products of Dehn twists on $S$.

The {\em conjugacy problem} asks for a decision procedure to 
determine whether or not $\gamma$ and $\phi$ are conjugate
in the mapping class group of $S$. That is, does there exist
an automorphism, $\varphi$, of $S$ such that 
$\gamma = \varphi^{-1} \phi \varphi$? 
\end{prob}

The {\em restricted conjugacy problem} is slightly easier in that 
it assumes that $\phi$ is in fact pseudo-Anosov. As this is
not a large restriction from now on we will assume that $\phi$
is pseudo-Anosov. 

A complete (albeit impractical) solution for the conjugacy problem 
has been given by Hemion~\cite{MR80f:57003}. The restricted case 
has also been solved by Mosher~\cite{MR89f:57016}.

\subsection{Notation}

In order to describe our solution of the restricted conjugacy problem, we
need to introduce some notation. Let $\gamma \co S\rightarrow S$ and
$\phi \co S\rightarrow S$ be pseudo-Anosov homeomorphisms with isometric
mapping tori. \snappea~will detect this and compute an isometry
$h\co M_\gamma \rightarrow M_\phi$.

Let $i_\gamma \co S \rightarrow M_\gamma$ and $i_\phi$ be the two inclusion
maps that realize $S$ as a fiber of the induced fiber structures 
$\mathcal{F}_\gamma$ and $\mathcal{F}_\phi$ on $M_\gamma$ and $M_\phi$.
Set $F_\gamma = i_\gamma(S)$ and define $F_\phi$ in a similar fashion.

Let $p_\gamma \co M_\gamma \rightarrow S^1$ be the map from 
$M_\gamma$ to the circle induced by $\mathcal{F}_\gamma$ 
and define $p_\phi$ similarly. 

Let $\sigma \in \Mod(M_\phi)$ denote a typical element of the isometry 
group of $M_\phi$. As $M_\phi$ is hyperbolic, $\Mod(M_\phi)$ is a
finite group and can be computed by \snappea.
Set $G_\sigma = \sigma h(F_\gamma)$. 

Finally, pick any $g \in \pi_1(M_\gamma)$ with the property that 
$(p_\gamma)_\ast g = 1 \in \mathbb{Z}$. We say that such a loop 
{\em represents the $S^1$-orientation}. 
If $(p_\phi \sigma h)_\ast g$ equals $1$ $(-1)$ then we say that $\sigma h$
{\em preserves (reverses) $S^1$-orientation}. 

\subsection{Retriangulation and the fundamental group}

Before solving the restricted conjugacy problem we
will need a pair of subroutines which determine the images of
elements of $\pi_1(M_\gamma)$ under the map $(p_\phi \sigma h)_\ast$.

First, we need an algorithm which decides whether 
$(p_\phi)_\ast \co \pi_1(G_\sigma) \rightarrow \mathbb{Z}$
has trivial image. The idea here is to keep track of
a set of generators for $\pi_1(S)$ under the maps
$i_\gamma$, $h$, $\sigma$ and $p_\phi$. Unfortunately, these homeomorphisms
do not all respect common simplicial structures on $M_\gamma$ and $M_\phi$.
To fix this problem one must find an appropriate set of generators 
in each retriangulation of $M_\gamma$. Then, once \snappea~finds a
geometric triangulation of $M_\gamma$, 
we can push the generators onto the one-skeleton
of the Ford domain. The isometry $\sigma h$ takes them to edge paths in the
one-skeleton of $M_\phi$ where we can reverse the process. 
Finally, $p_\phi$ projects the generators of $\pi_1(G_\sigma)$
to the circle where it is easy to check whether or not they are all
contractible.

Second, we will need an algorithm which decides whether $\sigma h$ 
preserves $S^1$-orientation. To do this, construct any loop 
$g \in \pi_1(M_\gamma)$ which represents the $S^1$-orientation. 
As in the previous algorithm take the image of $g$ under the 
map $(p_\phi \sigma h)_\ast$. Check that this image is the positive 
generator of $\pi_1(S^1)$.

The bookkeeping problem of keeping track of surface 
subgroups of $\pi_1(M^3)$ under retriangulation
does not yet have an implemented solution. 
We would be very interested in the work of any reader who is 
willing to write such a program. It should be remarked that
the subroutines above do a little more work than is 
strictly necessary. It would suffice to keep track of a
two-chain representing the fiber of $M_\gamma$. Again, 
this is a straightforward problem which does not yet have
an implemented solution.

\subsection{The algorithm}

Begin with two homeomorphisms $\gamma, \phi \co S\rightarrow S$.
If one of them fails to be pseudo-Anosov, then the software described
in \cite{pbexp} will detect this. Otherwise, construct the mapping tori
of $\gamma$ and $\phi$, $M_\gamma$ and $M_\phi$.
If \snappea~reports that $M_\gamma$ and $M_\phi$ are not isometric,
we conclude that $\gamma$ cannot be conjugate to 
$\phi$. This resolves the issue for a vast majority of possible
pairs of $\gamma$ and $\phi$.

Now, suppose \snappea~reports the two mapping tori are
isometric. We cannot yet conclude that the two 
automorphisms are conjugate. It may be that $\gamma$ and $\phi$
have homeomorphic mapping tori but are not conjugate because
they give rise to distinct fiber structures in the resulting 
three manifold. Also, it may happen that $\gamma$ and $\phi$ 
induce identical fiber structures while reversing $S^1$-orientation.
In this case we show that $\gamma$ is conjugate to $\phi^{-1}$. 

If $\pi_1(G_\sigma)$ has nontrivial image under $(p_\phi)_\ast$ then
clearly $G_\sigma$ is not isotopic to $F_\phi$.  We conclude that if
$\pi_1(G_\sigma)$ has nontrivial image in $\pi_1(S^1) = \mathbb{Z}$
for every $\sigma$ in $\Mod(M_\phi)$ then $\gamma$ is not conjugate to
$\phi$.

We claim that if there exists some $\sigma$ such that $\sigma h$ preserves 
$S^1$-orientation and $\pi_1(G_\sigma)$ is contained in the kernel 
of the natural projection, then $\gamma$ is conjugate to $\phi$,
which completes our solution of the restricted conjugacy problem.
The rest of this section is devoted to a proof of this claim.

\subsection{Isotoping $G_\sigma$}

Assume now that $\pi_1(G_\sigma)$ has trivial image in $(p_\phi)_\ast$
for some fixed $\sigma \in \Mod(M_\phi)$. At this point we need a weak 
form of Theorem~4 from \cite{ThurstonNorm}:

\begin{thmnp}
$G_\sigma$ is isotopic to a properly embedded surface that is
either a leaf of $\mathcal{F}_\phi$, or has only saddle singularities for 
the induced singular foliation of $G_\sigma$. The boundary 
component of the isotoped $G_\sigma$ is either a leaf of 
$\mathcal{F}_\phi | \partial M_\phi$ or is transverse to 
$\mathcal{F}_\phi | \partial M_\phi$.
\end{thmnp}

We use this as follows:

\begin{cor}
If $\pi_1(G_\sigma)$ is in the kernel of $(p_\phi)_\ast$
then $G_\sigma$ is isotopic to $F_\phi$.
\end{cor}

\begin{proof}
Suppose, to obtain a contradiction, that $G_\sigma$ is not
isotopic to $F_\phi$. By Thurston's theorem we may isotope
$G_\sigma$ so that the induced foliation has only saddle singularities.

Let $M_\mathbb{Z} = F_\phi \times \mathbb{R}$ 
be the infinite cyclic cover of $M$ coming from $\mathcal{F}_\phi$. 
By assumption we may lift $G_\sigma$ to $M_\mathbb{Z}$. 
Note that projection onto the second factor
$M_\mathbb{Z} \rightarrow \mathbb{R}$ 
gives a Morse function when restricted to $G_\sigma$.
However this Morse function must have a maximum.
The induced foliation of $G_\sigma$ has a singularity at this maximum,
but this singularity cannot possibly be a saddle singularity. This is
a contradiction.

\end{proof}

Thus we may isotope $\sigma h$ so as to obtain $G_\sigma = F_\phi$.
Cutting along $F_\gamma$ and $F_\phi$ we obtain a map $h' \co S \times I
\rightarrow S \times I$ that takes $S \times 0$ and $S \times 1$ to
$S \times 0$ and $S \times 1$, but not necessarily in that order. It may be 
that $\sigma h$ reverses the $S^1$-orientation.

It follows that $h'$, and hence $\sigma h$, is isotopic to a map that 
preserves fibers. (See, for example, Lemma 3.5 of~\cite{Waldhausen}.) 
Letting $h_0 = (\sigma h)|F_\gamma$ we find that either 
$\gamma = h_0^{-1} \phi h_0$ or $\gamma = h_0^{-1} \phi^{-1} 
h_0$ depending on whether $\sigma h$ preserves or reverses
$S^1$-orientation. This completes the proof of the claim and thus
shows the correctness of our algorithm, which we sum up in the following
theorem.

\begin{thmnp}\label{conj}
Let $\gamma$ and $\phi$ be pseudo-Anosov homeomorphisms of a once
punctured surface, $S$. The following is a procedure to decide whether the
two mappings are conjugate in the mapping class group of $S$, if~\snappea~
is allowed as a subroutine.

\begin{enumerate}
\item Apply Theorem~\ref{main} to obtain the mapping tori, $M_\gamma$ and
$M_\phi$.
\item Using~\jsnap~and~\snappea, determine whether or not $M_\gamma$ and
$M_\phi$ are isometric. If not, then $\gamma$ is not conjugate to $\phi$.
\item Using~\snappea, enumerate all isometries between $M_\gamma$ and
$M_\phi$.
\item Determine whether any of these are $S^1$-orientation and fibre
preserving. If none are, then $\gamma$ is not conjugate to $\phi$.
Otherwise, $\gamma$ is conjugate to $\phi$.
\end{enumerate}

\end{thmnp}

The complexity of this algorithm depends on the complexity
of \snappea, which we treat as a black box in this paper.

\appendix

\section{Generating input for \snappea}\label{format}


In order to generate input for \snappea, the output of \jmt~has
to be translated into \snappea's triangulation file format by a second
program (called \jsnap).
The purpose of this section is to discuss the intermediate format,
which allows for quick and easy generation of input for \snappea.
We illustrate this format with a familiar example from \cite{notes}.


The intermediate format admits two types of input lines, for tetrahedra
and for gluings. An input line defining a tetrahedron has the form
``{\tt T~$v_1~v_2~v_3~v_4$},'' where $v_1,\ldots,v_4$ are distinct labels of
the vertices. Tetrahedra are glued implicitly if they have three vertices
in common, and they can be glued explicitly by entering a line of the form
``{\tt G~$v_1~v_2~v_3~w_1~w_2~w_3$},'' where $v_1,v_2,v_3$ and $w_1,w_2,w_3$
are the labels of two faces of tetrahedra. In this gluing, the side $[v_1,v_2]$
is glued to the side $[w_1,w_2]$, the side $[v_2,v_3]$ is glued to the side
$[w_2,w_3]$, etc. Empty lines and comments beginning with `//' are also
allowed.

Although the four vertex labels of a tetrahedron have to be distinct, two
or more vertices of a tetrahedron may be identified after gluing. Distinct
vertex labels are only needed in order to uniquely specify the sides of
tetrahedra.

\begin{ex}[Figure-eight knot]\label{fig81}
Consider the familiar figure-eight knot (see \figref{fig8knot}). An ideal
triangulation of the complement of this knot can be expressed as a gluing
of two tetrahedra \cite[Chapter~1]{notes} (see \figref{fig8triang}).
\begin{figure}[tb]
\renewcommand{\epsfsize}[2]{0.4\textwidth}
\centerline{\epsfbox{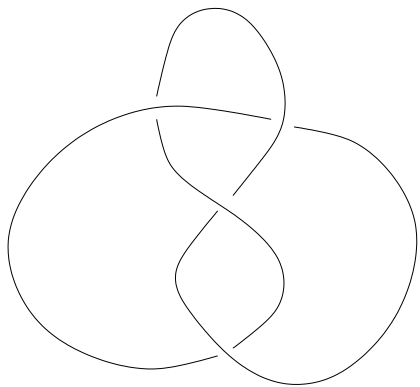}}
\caption{The figure-eight knot.}\label{fig8knot}
\end{figure}

\begin{figure}[tb]
\renewcommand{\epsfsize}[2]{0.7\textwidth}
\psfrag{a}{$a$}
\psfrag{b}{$b$}
\psfrag{c}{$c$}
\psfrag{d}{$d$}
\psfrag{e}{$e$}
\centerline{\epsfbox{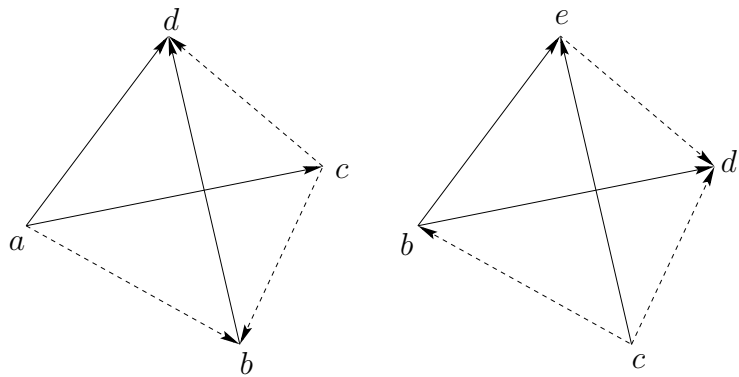}}
\caption{The figure-eight knot complement, expressed as a gluing of two
ideal tetrahedra. Triangles are glued such that the line style and the
direction of arrows are matched \cite[Chapter~1]{notes}.}
\label{fig8triang}
\end{figure}

We encode the gluing as follows.
\begin{verbatim}
T a b c d
T b c d e
G b e d a c d
G c b e a b d
G c e d a c b
\end{verbatim}%

Feeding the above five lines into \jsnap~yields an encoding of the
triangulation in \snappea's file format. According to \snappea, the
fundamental group of the resulting manifold $M$ has the following
presentation, which we modify by a sequence of Tietze transformations
\cite[Chapter~II.2]{ls}.
\begin{align}
\pi_1(M)~&\cong~<x,~y~|~\bar y\bar x\bar x\bar x\bar yxyyx=1>\label{fig8snap}\\
&\cong~<x,~y~,~z~|~\bar y\bar x\bar x\bar x\bar yxyyx=1,~y=z\bar x>\notag\\
&\cong~<x,~z~|~x\bar z\bar x\bar x\bar zxz\bar xz=1>\notag\\
&\cong~<x,~z~|~\bar xzx\bar z\bar x\bar x\bar zxz=1>\label{fig8gp}
\end{align}
\presref{fig8gp} is the presentation of the fundamental group of
the complement of the figure-eight knot given in \cite[Example~3.8]{knots}.
\end{ex}

While \snappea's format is an extremely efficient representation of
triangulations, understanding it requires some effort. The program
\jsnap~acts as an interface between \snappea~and the human user.
If you can draw or visualize a triangulation of a $3$-manifold,
then you can also enter it into \jsnap.

\section{Sample Computations}\label{exsec}

\begin{figure}[tb]
\renewcommand{\epsfsize}[2]{0.8\textwidth}
\psfrag{a0}{$a_0$}
\psfrag{b0}{$b_0$}
\psfrag{c0}{$c_0$}
\psfrag{d0}{$d_0$}
\psfrag{a1}{$a_1$}
\psfrag{b1}{$b_1$}
\psfrag{c1}{$c_1$}
\psfrag{d1}{$d_1$}
\centerline{\epsfbox{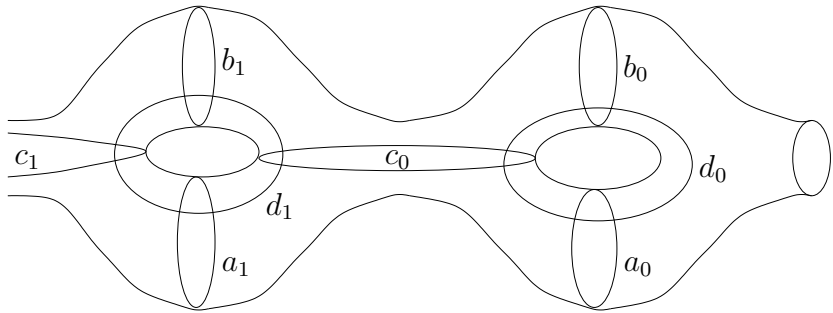}}
\caption{A set of generators of the mapping class group.}
\label{stdtwist}
\end{figure}

We present some sample computations that illustrate the power of the software
discussed here. We define surface automorphisms as compositions of Dehn twists
with respect to the set of curves shown in \figref{stdtwist}.
We equip the surface with a normal vector field. When twisting with respect
to a curve $c$, we turn right whenever we hit $c$ (resp.\ left for inverse
twists). The set of these Dehn twists generates the mapping class group
\cite{lickorish}.

We first present an example for which we can easily verify correctness.

\begin{ex}[Figure-eight knot revisited]\label{fig82}
Let $S$ be a punctured torus, and let $\phi\co S\rightarrow S$ be given by
$$\phi=D^{-1}_{a_0}D^{\null}_{d_0}.$$
The complement of the figure-eight knot (see \figref{fig8knot}) is
homeomorphic to the mapping torus of $\phi$ (we will verify this soon).
Given this composition of Dehn twists, the train track software described in 
\cite{pbexp} determines that $\phi$
is pseudo-Anosov with growth rate $\lambda\approx 2.61803399,$ and
\snappea~determines that the mapping torus $M_\phi$ is a hyperbolic
$3$-manifold of volume $V\approx 2.02988321$ with one torus cusp.

Moreover, the train track software computes the following topological
representative $f\co G\rightarrow G$ of $\phi$, where $G$ is a graph with one
vertex and two edges.
\begin{eqnarray*}
f(a)&=&ba\\
f(b)&=&bba\\
\sigma&=&a\bar b\bar ab\\
\end{eqnarray*}
This yields the following presentation of the fundamental group of the
mapping torus $M_\phi$ of $\phi$.
$$ \pi_1(M_\phi)~\cong~<a,~b,~t~|~\bar tat=ba,~\bar tbt=bba> $$
A few Tietze transformations
show that $\pi_1(M_\phi)$ is the fundamental group of the figure-eight
knot complement.
\begin{align*}
\pi_1(M_\phi)~&\cong~<a,~b,~t~|~\bar tat=ba,~\bar tbt=bba>\\
&\cong~<a,~t~|~t\bar a\bar ata\bar t\bar a\bar ta = 1>\\
&\cong~<a,~c~|~\bar aca\bar c\bar a\bar a\bar cac = 1>
\end{align*}
As in \exref{fig81}, we have obtained the presentation given in
\cite[Example~3.8]{knots}.
Finally, the presentation of $\pi_1(M_\phi)$ computed by \snappea~is
$$\pi_1(M_\phi)~\cong~<x,~y~|~\bar x\bar y\bar y\bar y\bar xyxxy = 1>,$$
which agrees with \presref{fig8snap} in \exref{fig81}.

Alternatively, we can run \snappea's isometry checker on \exref{fig81} and
\exref{fig82} in order to see that we get the same hyperbolic $3$-manifold
in both examples.
\end{ex}
\begin{ex}[Genus 3]
Let $S$ be a surface of genus $3$ with one puncture, and let
$\phi\co S\rightarrow S$ be given by
$$\phi=D^{\null}_{d_0}D^{\null}_{c_0}D^{\null}_{d_1}D^{\null}_{c_1}D^{\null}_{d_2}D^{-1}_{a_2}.$$
The train track software identifies $\phi$ as a pseudo-Anosov
homeomorphism with growth rate $\lambda\approx 2.04249053$, and \snappea~
determines that the mapping torus $M_\phi$ is a hyperbolic $3$-manifold of
volume $V\approx 4.93524268$ with one torus cusp.
\end{ex}

These applications only show a small part of all the possibilities. The train
track software computes a plethora of information about surface homeomorphisms
\cite{pbexp}, and \snappea~allows for a detailed analysis of (hyperbolic)
$3$-manifolds. We believe that the combination of the two packages may
become a valuable tool for topologists.

\bibliographystyle{alpha}
\bibliography{$GLOBAL/all}
\par

{\noindent \sc Peter Brinkmann\\
Department of Mathematics\\
273 Altgeld Hall\\
1409 W.\ Green Street\\
Urbana, IL 61801, USA\\
E-mail: brinkman@math.uiuc.edu\\
\\}
{\noindent \sc Saul Schleimer\\
Department of Mathematics\\
UC Berkeley\\
Berkeley, CA 94720, USA\\
E-mail: saul@math.berkeley.edu\\
\\}

\end{document}